\documentclass[11pt]{article}
\usepackage{amssymb}
\usepackage{amsmath}
\usepackage{amsthm}
\usepackage{latexsym}
\usepackage{amsfonts}
\usepackage{graphicx}
\usepackage{graphics}

\newcommand{\dis}{\displaystyle}
\theoremstyle{plain}
\newtheorem{thm}{Theorem}[section]   % Αρίθμηση συνεχόμενη (όχι κατά Θεώρημα, Λήμμα κ.λπ.)
\newtheorem{prop}[thm]{Proposition}

\newtheorem{cor}[thm]{Corollary}
\newtheorem{Def}[thm]{Definition}%[section]

\theoremstyle{definition}
\newtheorem{rem}[thm]{Remark}

\newtheorem*{Proof}{Proof}

\newcommand{\cU}{\mathcal{U}}

\newcommand{\el}{\ell}

\newcommand{\ra}{\;\rightarrow\;}

\newcommand{\ga}{\gamma }

\newcommand{\de}{\delta }
\newcommand{\OO} {{\varOmega}}
\newcommand{\De} {{\varDelta}}
\newcommand{\e}{\varepsilon }

\newcommand{\zi}{\zeta }

\newcommand{\vthi}{\vartheta }

\newcommand{\la}{\lambda }
\newcommand{\mi}{\mu }

\newcommand{\oo}{\omega}
\newcommand{\C}{\mathbb{C}}

\newcommand{\N}{\mathbb{N}}
\newcommand{\Q}{\mathbb{Q}}
\newcommand{\F}{\mathbb{F}}

\newcommand{\ssum}{\sum\limits}

\newcommand{\oD}{\overline{D}}

\newcommand{\oV}{\overline{V}}
\newcommand{\oS}{\overline{S}}
\newcommand{\oB}{\overline{B}}
\newcommand{\oT}{\overline{T}}
\newcommand{\oL}{\overline{L}}
\newcommand{\oO}{\overline{\varOmega}}

\newcommand{\tp}{\widetilde{p}}
\newcommand{\tq}{\widetilde{q}}

\newcommand{\tA}{\widetilde{A}}
\newcommand{\tB}{\widetilde{B}}

\newcommand{\cu}{{\cal{U}}}

\newcommand{\cD}{{\cal{D}}}

\newcommand{\ld}{\ldots}

\newcommand{\sm}{\smallsetminus}

\newcommand{\qb}{$\quad\blacksquare$}

 \newcommand{\Arg}{\mbox{Arg}}

 \renewcommand{\Re}{\mbox{Re}}

  \newcommand{\dsup}{\displaystyle\sup}
 
 \newcommand{\dist}{\mbox{dist}}
\begin{document}
\title{\bf Universal Pad\'{e} approximants and their behaviour on the boundary}
\author{Ilias Zadik}
\date{}
\maketitle
\begin{abstract}
There are several kinds of universal Taylor series. In one such kind the universal approximation is required at every boundary point of the domain of definition $\OO$ of the universal function $f$. In another kind the universal approximation is not required at any point of $\partial\OO$ but in this case the universal function $f$ can be taken smooth on $\oO$ and, moreover, it can be approximated by it's Taylor partial sums on every compact subset of $\oO$. Similar generic phenomena hold when the partial sums of the Taylor expansion of the universal function are replaced by some Pad\'{e} approximants of it.

In the present paper we show that in the case of Pad\'{e} approximants, if $\OO$ is an open set and $S,T$ are two subsets of $\partial\OO$ that satisfy some conditions, then there exists a universal function $f\in H(\OO)$ which is smooth on $\OO\cup S$ and has some Pad\'{e} approximants that approximate $f$ on each compact subset of $\OO\cup S$ and simultaneously obtain universal approximation on each compact subset of $(\C\sm\oO)\cup T$. A sufficient condition for the above to happen is $\oS\cap\oT=\emptyset$, while a necessary and sufficient condition is not known.
\end{abstract}
\section{Introduction}\label{sec1}
\noindent

Universality is a generic phenomenon of chaotic behavior of approximation \cite{2}. One can say that a universal series of functions is a series of functions that are so extremely divergent that subsequences of them approximate lots of different functions on some sets.

A basic example of this phenomenon is the case of universal Taylor series, where the partial sums of the Taylor expansion of a holomorphic function $f$ approximate many functions on some compact sets outside the domain of definition $\OO$ of $f$, with respect to the Euclidean distance. The history started with Fekete \cite{15} before 1914, who introduced the first universal Taylor series on the real line. In the early 70'\,s Luh \cite{9} and independently Chui and Parnes \cite{3} proved the existence of a Taylor series around zero with radius of convergence equal to 1 whose partial sums approximate uniformly any polynomial function on any compact set $K\subset\C\mid\oD=\{z\in\C\mid|z|>1\}$ with connected complement. In \cite{11}, Nestoridis showed that there exists a universal Taylor series around zero with radius of convergence equal to 1, where the approximation we described above is valid also on compact sets $K$ that can touch the boundary $\partial D=\{z\in\C\mid|z|=1\}$. Moreover, he proved in \cite{12} that generically every function $f\in H(\OO)$, where $\OO$ is a simply connected domain, has universal Taylor series with respect to any center $\zi\in\OO$.

These series have universal approximating properties on every compact set $K$ in $\C\sm\OO$ with connected complement.

It is also proved in \cite{7} that there exists a function $f\in H(\OO)$ with the above properties which is also smooth on $\partial\OO$, provided that the compact sets $K$ in $\C\sm\OO$ where the universal approximation takes place do not touch the boundary.

Recently, the partial sums of the Taylor expansion of a function, which are polynomials, have been replaced by a set of Pad\'{e} approximants of $f$, which are rational functions. Therefore, they can take the value $\infty$ as well. In accordance with universal Taylor series, several results were proved in this new direction(\cite{4}, \cite{5}, \cite{13}, \cite{14}). In these new results it is shown that the Pad\'{e} approximants  can approximate any rational function uniformly with respect to the chordal metric on any compact set contained in $\C\sm\oO$ \cite{14}, or in $\C\sm\OO$ \cite{13}. Moreover, it was proved in \cite{6} that generically every function $f\in H(\OO)$ can be approximated on any compact set in $\OO$ by it's Pad\'{e} approximants with respect to the Euclidean metric this time.

Thus, the above results assure us that generically, every function in $H(\OO)$ can be approximated by some of it's Pad\'{e} approximants and, simultaneously, the same Pad\'{e} approximants realize universal Pad\'{e} approximation of any rational function outside of $\oO$.

The question that naturally arises is what can be achieved on the boundary of $\OO$. In \cite{13}, the universal approximation was extended to any compact subset of $\partial\OO$ while in \cite{14} the approximation of $f$ by its Pad\'{e} approximants is valid not only on $\OO$ but on $\partial\OO$ as well.

In the present paper we show that $\partial\OO$ can be splitted in two disjoint parts $S$ and $T$, which must satisfy some conditions, so that the approximation towards $f$ is valid on compact subsets of $\OO\cup S$ and the universal approximation of any rational function is valid on compact subsets of $(\C\sm\oO)\cup T$. Our result is a general one and yields as corollaries all previous results. It is also generic in a closed subspace of $T^\infty(\OO,(L_n)_{n\in\N})$ where $L_n$, $n=1,2,\ld$ are compact subsets of $\oO$ with some specific properties and $T^\infty(\OO,(L_n)_{n\in\N})$ is the set of all functions $f\in H(\OO)$ where $f^{(\el)}
\mid_{(L_n\cap\OO)}$ is uniformly continuous for every $\el\in\N_0$, $n\in\N$. The topology we consider on $T^\infty(\OO,(L_n)_{n\in\N})$ is endowed by the seminorms $\{\dis\sup_{z\in(L_n\cap\OO)}|f^{(\el)}(z)|\mid\el\in\N_0,n\in\N\}$ and therefore $T^\infty(\OO_n(L_n)_{n\in\N})$ becomes a Frechet space.

The closed subspace of $T^\infty(\OO)=T^\infty(\OO,(L_n)_{n\in\N})$ where our result is generic, consists of the closure in $T^\infty(\OO,(L_n)_{n\in\N})$ of all rational functions with poles off $\dis\bigcup_{n\in\N}L_n$.

Furthermore, we consider compact sets $K_m\subset\C$, $m=1,2,\ld$, so that $K_m\cap L_n=\emptyset$ for all $m,n\in\N$ and we obtain universal Pad\'{e} approximation on $K_m$, $m=1,2,\ld$ with respect to the chordal metric. For special choices of the sequences $(L_n)_{n\in\N}$ and $(K_m)_{m\in\N}$ we cover all previously obtained generic results.

Finally, in the last section of the present paper we show that if the boundary $\partial\OO$ of an open set $\OO\subseteq\C$ contains two parts $S,T$ with $\oS\cap\oT=\emptyset$, then there always exists a function $f\in H(\OO)$ such that some Pad\'{e} approximants of $f$ approximate any rational function on the compact subsets of $(\C\sm\oO)\cup T$ and simultaneously approximate $f$ on the compact subsets of $\OO\cup S$, where $f$ is smooth.
\section{Preliminaries}\label{sec2}
\noindent

If we consider the one-point compactification
$\C\cup\{\infty\}=\widetilde{\C}$ of $\C$, then a well known metric
is the chordal metric $\chi$ on $\C\cup\{\infty\}$, where
\[
\chi(a,b)=\frac{|a-b|}{\sqrt{1+|a|^2}\sqrt{1+|b|^2}}, \ \ \text{for} \ \ a,b\in\C
\]
and $\chi(a,\infty)=\dfrac{1}{\sqrt{1+|a|^2}}$ for $a\in\C$, and
$\chi(\infty,\infty)=0$;
\begin{prop}\label{prop2.1}
Let $K\subset\C$ be a compact set and $q=\dfrac{A}{B}$ a rational
function, where the polynomials $A,B$ do not have a common zero in
$\C$.\ Then there is a sequence $q_j=\dfrac{A_j}{B_j}$, $j=1,2,\ld$
where the polynomials $A_j$ and $B_j$ have coefficients in $\Q+i\Q$
and do not have any common zero in $\C$ for all $j$, such that
$\dis\sup_{z\in K}\chi(q_j(z),q(z))\ra0$ as $j\ra+\infty$.
\end{prop}

The above proposition is well known.\ See \cite{13}.

Let $\zi\in\C$ be fixed and
\[
f=\sum^\infty_{n=0}a_n(z-\zi)^n
\]
be a formal power series $(a_n=a_n(f,\zi))$.\ This power series is often the Taylor development of a holomorphic function $f$ in a neighborhood of $\zi$.\ Let $p$ and $q$ be two non negative integers.\ The Pad\'{e} approximant $[f;p/q]_\zi(z)$ is defined to be a rational function $\phi$ regular at $\zi$ whose Taylor development with center $\zi$,
\[
\phi(z)=\sum^\infty_{n=0}b_n(z-\zi)^n,
\]
satisfies $b_n=a_n$ for all $0\le n\le p+q$ and $\phi(z)=A(z)/B(z)$, where the polynomials $A$ and $B$ satisfy
\[
\deg A\le p, \ \ \deg B\le q \ \ \text{and} \ \ B(\zi)\neq0.
\]
It is not always true that such a rational function $\phi$ exists.\ And if it exists it is not always unique.\ For $q=0$, we always have such a unique $\phi$ which is
\[
[f;p/q]_\zi(z)=\sum^p_{n=0}a_n(z-\zi)^n.
\]
For $q\ge1$ the necessary and sufficient condition for existence and
uniqueness is that the following $q\times q$ Hankel determinant is
non-zero (\cite{1})
\[
\left|\begin{array}{lllll}
        a_{p-q+1} & a_{p-q+2} & \cdot & \cdot & a_p \\
        a_{p-q+2} & a_{p-q+3} & \cdot & \cdot & a_{p+1} \\
        \cdot & \cdot & \cdot & \cdot & \cdot \\
        \cdot & \cdot & \cdot & \cdot & \cdot \\
        a_p & a_{p+1} & \cdot & \cdot & a_{p+q-1}
      \end{array}\right|\neq0,
\]
where $a_i=0$ for $i<0$.\ If this is satisfied we write $f\in\cD_{p,q}(\zi$).\ For $f\in\cD_{p,q}(\zi)$ the Pad\'{e} approximant
\[
[f;p/q]_\zi(z)=\frac{A(f,\zi)(z)}{B(f,\zi)(z)}
\]
is given by the following Jacobi formula
\[
A(f,\zi)(z)=\left|\begin{array}{lllll}
                    (z-\zi)^qS_{p-q}(f,\zi)(z) & (z-\zi)^{q-1}S_{p-q+1}(f,\zi)(z) & \cdot & \cdot & S_p(f,\zi)(z) \\
                    a_{p-q+1} & a_{p-q+1} & \cdot & \cdot & a_{p+1} \\
                    \cdot & \cdot & \cdot & \cdot & \cdot \\
                     \cdot & \cdot & \cdot & \cdot & \cdot \\
                    a_p & a_{p+1} & \cdot & \cdot & a_{p+q}
                  \end{array}\right|,
\]
\[
B(f,\zi)(z)=\left|\begin{array}{ccccc}
                    (z-\zi)^q & (z-\zi)^{q-1} & \cdot & \cdot & 1 \\
                    a_{p-q+1} & a_{p-q+1} & \cdot & \cdot & a_{p+1} \\
                    \cdot & \cdot & \cdot & \cdot & \cdot \\
                    \cdot & \cdot & \cdot & \cdot & \cdot \\
                     a_p & a_{p+1} & \cdot & \cdot & a_{p+q}
                  \end{array}\right|,
\]
with (see \cite{1})
\[
S_k(f,\zi)(z)=\left\{\begin{array}{ccc}
                       \ssum^k_{\nu=0}a_\nu(z-\zi)^\nu, & \text{if} & k\ge0 \\
                       0, & \text{if} & k<0.
                     \end{array}\right.
\]
If $A(f,\zi)(z)$ and $B(f,\zi)(z)$ are given by the previous Jacobi formula and do not have a common zero in a set $K$ we write $f\in E_{p,q,\zi}(K)$.\ Equivalently
\[
|A(f,\zi)(z)|^2+|B(f,\zi)(z)|^2\neq0
\]
for all $z\in K$.\ For $K$ compact this is equivalent to the existence of a $\de>0$ such that
\[
|A(f,\zi)(z)|^2+|B(f,\zi)(z)|^2>\de
\]
for all $z\in K$.\ We will also use the following ({\cite{1}
Th.\ 1.4.4 page 30).
\begin{prop}\label{prop2.2}
Let $\phi(z)=\dfrac{A(z)}{B(z)}$ be a rational function, where the
polynomials $A$ and $B$ do not have any common zero in $\C$.\ Let
$\deg A(z)=k$ and $\deg B(z)=\la$.\ Then for every $\zi\in\C$ such
that $B(\zi)\neq0$ we have the following:
\[
\phi\in\cD_{k,\la}(\zi),
\]
\[
\phi\in\cD_{p,\la}(\zi) \ \ \text{for all} \ \ p>k,
\]
\[
\phi\in\cD_{k,q}(\zi) \ \ \text{for all} \ \ q>\la.
\]
In all these cases $\phi$ coincides with its corresponding Pad\'{e}
approximant, that is,\linebreak $[\phi;k/\la]_\zi(z)\equiv\phi(z)$
and $[\phi;p/\la]_\zi(z)\equiv\phi(z)$ for all $p>k$ and
$[\phi;k/q]_\zi(z)\equiv\phi(z)$ for $q>\la$.
\end{prop}

We shall also make use of Lemma 1 from \cite{10} which states the following.
\begin{prop}\label{prop2.3}
Let $\OO\subset\C$ be an open set. We suppose that the number of components of $\OO$ is locally finite. Then there exists a sequence $K_m\subset\C\sm\OO$ $m=1,2,\ld$ of compact sets such that $K_m\cap\OO=\emptyset$ with $K^c_m$ connected such that for every $K\subset\C\sm\OO$ compact set with $K^c$ connected there exists $m=1,2,\ld$ with $K\subset K_m$.
\end{prop}

Moreover, a slight modification at the proof of the above proposition that appears in \cite{10} gives the following result.
\begin{prop}\label{prop2.4}
Let $\OO\subseteq\C$ be an open set. We suppose that the number of components of $\OO$ is locally finite. Then there exists a sequence $K_m\subset\C\sm\oO$, $m=1,2,\ld$ of compact sets with $K^c_m$ connected such that for every compact set $K\subset\C\sm\oO$ with $K^c$ connected there exists $m=1,2,\ld$ with $K\subset K_m$.
\end{prop}
\begin{Def}\label{Def2.5}
Let $\OO\subseteq\C$. We say that $(L_n)_{n\in\N}$ is an exhausting sequence of compact sets in $\OO$ if $\OO=\dis\bigcup_{n\in\N}L_n$ and for every $n\in\N$, it holds that $L_n\subset L^0_{n+1}$, where the interior is taken with respect to the relative topology of $\OO$.
\end{Def}
\begin{prop}\label{prop2.6}
If $(L_n)_{n\in\N}$ is an exhausting sequence of compact sets in $\OO\subseteq\C$, then if $J\subset\OO$ is a compact set, there exists $N\in\N$ such that $J\subset L_N$.
\end{prop}
\begin{Proof}
It is $J\subset\OO=\dis\bigcup_{n\in\N}L_n\subset\dis\bigcup_{n\ge2}L^0_n$ and therefore there exists $N\in\N$ such that $J\subset\dis\bigcup^N_{n=2}L^0_n$, because $\zi$ is compact. But $(L_n)_{n\in\N}$ is an increasing sequence of compact sets and therefore it holds $\zi\subset\dis\bigcup^N_{n=2}L^0_n=L^0_N\subset L_N$, as we wanted.\qb
\end{Proof}
\begin{prop}\label{prop2.7}
If $\OO\subseteq\C$ is locally compact then there exists an exhausting sequence of compact sets in $\OO$.
\end{prop}
\begin{Proof}
Let $\tB$ be a countable basis for $\OO$ and consider $V=\{B\in\tB\mid\oB$ is compact$\}$.

Then $V$ is a basis since for any open $\cU$ and $x\in\cU$ we choose a compact neighborhood $N$ of $x$ inside $U$ and because $\tB$ is a basis there exists $B\in\tB$ such that $x\in B\subset N$; this implies $\oB\subset N$ which yields that $B$ is compact and we found a $B\in V$ such that $x\in B\subset\cU$.

Therefore, we may assume that we have a countable basis $V=\{V_n\mid n\in\N\}$ where $\oV_n$ is compact for every $n\in\N$. We define now the exhausting sequence of compact sets as follows.

We put $L_1=\oV_1$. Since $V$ covers the compact set $L_1$ there exist $i_1\in\N$ such that $L_1\subset V_1\cup V_2\cup\cdots\cup V_{i_1}$.

Denoting $D_1$ the right hand side of the inclusion above, we put $L_2=\oV_1\cup\oV_2\cup\cdots\cup\oV_{i_2}=\oD_1$.

This set is compact as a finite union of compacts.

Since $D_1\subset L_2$ and $D_1$ is open we get $D_1\subset L^0_2$, which gives $L_1\subset L^0_2$.

Next we choose $i_2>i_1$ such that $L_2\subset V_1\cup V_2\cup\cdots\cup V_{i_2}$, denote $D_2=V_1\cup V_2\cup\cdots \cup V_{i_2}$ and put $L_3=\oD_2=\oV_1\cup\oV_2\cup\cdots\cup\oV_{i_2}$. As before $L_3$ is compact and $L_2\subset L^0_3$.

Continuing this process we construct the sequence $(L_n)_{n\in\N}$ which clearly covers $\OO$ and for every $n\in\N$ it holds $L_n\subset L^0_{n+1}$. \qb
\end{Proof}
\begin{prop}\label{prop2.8}
Let $\OO$ be an open set and $K\subseteq\C$ a compact subset of $\oO$. Then there exists a compact subset $K'$ of $\oO$ such that\vspace*{-0.2cm}
\begin{enumerate}
\item[(i)] $K\cap\partial\OO=K'\cap\partial\OO$\vspace*{-0.2cm}
\item[(ii)] $K\subseteq K'$ and\vspace*{-0.2cm}
\item[(iii)] Every connected component of $\{\infty\}\cup(\C\sm K')$ contains a point from $\{\infty\}\cup(\C\sm\oO)$.
\end{enumerate}
\end{prop}
\begin{Proof}
It is enough to show that if $(L_n)_{n\in I}$, $I\subseteq N$ are the bounded connected components of $\{\infty\}\cup(\C\sm K)$ that are fully contained in $\OO$, then the set $K'=\Big(\dis\bigcup_{n\in I}L_n\Big)\cup K$ is compact.

One can see that for every $n\in\N$, $\partial L_n\subset K$. This fact together with the fact that $K$ is bounded, yields that the union $\dis\bigcup_{n\in I}L_n$ must also be bounded. Therefore, it is enough to show that $K'$ is closed. Let $(x_n)_{n\in\N}\subset K'$ with $x_n\ra x$.

Since $K'\subset\oO$, $x\in\oO$. Assume now that $x\notin K'$. Then of course $x\notin K$ and since $K$ is compact, eventually $x_n\in\dis\bigcup_{N\in I}L_N$. Moreover, since $\oL_n=L_n\cup\partial L_n\subset L_n\cup K\subset K'$, it is also true that $(x_n)_{n\in\N}$ does not have a whole subsequence inside one $L_N$, $N\in I$. Therefore $\dist(x,L_{m_n})\ra0$, as $n\ra+\infty$, for a subsequence $(L_{m_n})_{n\in\N}$ of $(L_n)_{n\in\N}$. But this means $\dist(x,\partial L_{m_n})\ra0$ as $n\ra+\infty$. One can prove the above fact by observing that if $\la\in  L_N$ for some $N\in\N$ then there exist $\la'\in\partial L_N$ with $|\la'-x|\le|\la-x|$, because the line segment $[\la,x]$ is a connected set and therefore the equality
\[
[\la,x]=(([\la,x])\cap L^0_N)\cup([\la,x]\cap\partial L_N)\cup([\la,x]\cap(\C\sm L_N))
\]
yields $[\la,x]\cap\partial L_N\neq\emptyset$.

Now, the fact that $\dist(x,\partial L_{w_n})\ra0$, as $n\ra+\infty$ combined with $\partial L_{m_n}\subset K$ yields $\dist(x,K)=0$, which is a contradiction since $K$ is closed and $x\notin K$. \qb
\end{Proof}
\begin{cor}\label{cor2.9}
Let $\OO$ be an open set and $(L_n)_{n\in\N}$ be an increasing sequence of compact sets in $\oO$ such that \vspace*{-0.2cm}
\begin{enumerate}
\item[(i)] For every $n\in\N$, $\overline{(L_n\cap\OO)}=L_n$.\vspace*{-0.2cm}
\item[(ii)] For every compact set $J$ in $\OO$ there exists $N\in\N$ such that $J\subset L_N$.
\end{enumerate}\vspace*{-0.2cm}

Then, there exists an increasing sequence $(L'_n)_{n\in\N}$ of compact sets in $\oO$ such that\vspace*{-0.2cm}
\begin{enumerate}
\item[1)] $(L'_n)_{n\in\N}$ satisfy conditions (i), (ii).\vspace*{-0.2cm}
\item[2)] For every $n\in\N$ it holds $L_n\subseteq L'_n$, $L_n\cap\partial\OO=L'_n\cap\partial\OO$ and every connected component of $\{\infty\}\cup(\C\sm L'_n)$ contains a point from $\{\infty\}\cup(\C\sm\oO)$.
\end{enumerate}
\end{cor}
\begin{Proof}
From Proposition \ref{prop2.8} for every $L_n$, $n\in\N$, we can find a compact set $L'_n$ such that $L_n\subseteq L'_n$, $L_n\cap\partial\OO=L'_n\cap\partial\OO$ and every connected component of $\{\infty\}\cup(\C\sm L'_n)$ contains a point from $\{\infty\}\cup(\C\sm\oO)$.

One can check that the sequence $(L'_n)_{n\in\N}$ satisfy the conditions we want. \qb
\end{Proof}
\section{Main Results}\label{sec3}
\noindent

Let $\OO$ be an open set, $\OO\subseteq\C$ and $L_n$, $n\in\N$, be an increasing sequence of compact subsets of $\oO$. For this sequence we also assume firstly that for every $n\in\N$ it holds $\overline{(L_n,\cap\OO)}=L_n$, secondly that each connected component of $\{\infty\}\cup(\C\sm L_n)$ contains a connected component of $\{\infty\}\cup(\C\sm\oO)$ and finally that every compact set $J$ inside of $\OO$ is contained in some $L_m$, for a natural number $m$.

Let $T^\infty(\OO)=T^\infty(\OO,(L_n)_{n\in\N})$ be the space of all analytic functions $f\in H(\OO)$, such that, for every derivative $f^{(\el)}$, $\el\in\N_0$, and every $L_n$, $n\in\N$ the restriction $f^{(\el)}\mid_{(L_n\cap\OO)}$ is uniformly continuous and therefore it extends continuously on $\overline{(L_n\cap\OO)}=L_n$.

We endow this space with the seminorms $\dsup_{z\in L_n}|f^{(\el)}(z)|$, $\el\in\N_0$, $n\in\N$. Then, it becomes a Frechet space, containing the rational functions with poles off the union $\dis\bigcup_{n\in\N}L_n$. Consider now, $Y^\infty(\OO)$ the closure of the set of all the rational functions with poles off $\dis\bigcup_{n\in\N}L_n$. As a closed subset of a complete space, $Y^\infty(\OO)$ is also a complete space and therefore Baire's theorems is at our disposal.

Based on the above facts we are now ready for our first theorem.
\begin{thm}\label{thm3.1}
Let $F\subset\N\times\N$ be a set that contains a sequence $(\tp_n,\tq_n)$, $n=1,2,\ld$ such that $\tp_n\ra+\infty$ and $\tq_n\ra+\infty$ and let $\OO\subseteq\C$ be an open set. Let also $K\subseteq\C$ be a compact set such that $K\cap L_n=\emptyset$ for every $n\in\N$ and let $m\in\N$ be a fixed natural number.

Then there exists $f\in Y^\infty(\OO)$ such that: for every rational function $h$ there exists a sequence $(p_n,q_n)\in F$ $(n=1,2,\ld)$ with the following properties:\vspace*{-0.2cm}
\begin{enumerate}
\item[(i)] $f\in D_{p_n,q_n}(\zi)\cap E_{p_n,q_n,\zi}(L_m\cup K)$, for every $\zi\in L_m$, $n\in\N$.\vspace*{-0.2cm}
\item[(ii)] For every $\el\in\N$, $\dsup_{\zi\in L_m}
\dsup_{z\in L_m}\big|\big[f;p_n/q_n\big]_\zi^{(\el)}(z)-f^{(\el)}(z)\big|\ra0$, as $n\ra+\infty$.\vspace*{-0.2cm}
\item[(iii)] $\dsup_{\zi\in L_m}\dsup_{z\in K}\chi([f;p_n/q_n]_\zi(z),\,h(z))\ra0$, as $n\ra+\infty$.
\end{enumerate}\vspace*{-0.2cm}

The set of such functions $f\in Y^\infty(\OO)$ is dense and $G_\de$ in $Y^\infty(\OO)$.
\end{thm}
\begin{Proof}
Let $(f_j)_{j\in\N}$ be an enumeration of the rational functions having coefficients of the numerator and denominator from $\Q+i\Q$.

We denote $\cu^K(m)$ the set of all functions in $Y^\infty(\OO)$ that satisfy the properties (i), (ii) and (iii).

We will prove that $\cu^K(m)$ is a $G_\de$-dense subset of $Y^\infty(\OO)$ and therefore $\cu^K(m)\neq\emptyset$.

For $j,s\in\N$ and $(p,q)\in F$ we define:
\begin{align*}
E(j,p,q,s)=\big\{&f\in Y^\infty(\OO)\mid f\in D_{p,q}(\zi)\cap E_{p,q,\zi}(K) \ \ \text{for all} \ \ \zi\in L_m\\
&\text{and} \ \ \dsup_{\zi\in L_m}\dsup_{z\in K}\chi([f;p/q]_\zi(z),\,f_j(z))<1/s\big\}
\end{align*}
and,
\begin{align*}
T(p,q,s)=\big\{&f\in Y^\infty(\OO)\mid f\in D_{p,q}(\zi)\cap E_{p,q,\zi}(L_m) \ \ \text{for all} \ \ \zi\in L_m \\
&\text{and} \ \ \dsup_{\zi\in L_m}\dsup_{z\in L_m}\big|\big[f;p/q\big]^{(\el)}_\zi(z)-f^{(\el)}(z)\big|<1/s \ \ \text{for} \ \ \el=0,1,\ld,s\big\}.
\end{align*}
Proposition \ref{prop2.1} and the definition of $Y^\infty(\OO)$ easily implies that
\[
U^K(m)=\bigcap^{+\infty}_{j,s=1}\bigcup_{(p,q)\in F}(E(j,p,q,s)\cap T(p,q,s)).
\]
To prove that $U^K(m)$ is a $G_\de$-dense in the $Y^\infty(\OO)$-topology, it is enough to prove that for every $j,s=1,2,\ld$ and $(p,q)\in F$ the sets $E(j,p,q,s)$, $T(p,q,s)$ are open in $Y^\infty(\OO)$ and that for every $j$ and $s$ from $\N^\ast$, the set $\dis\bigcup_{(p,q)\in F}(E(j,p,q,s)\cap T(p,q,s))$ is dense in $Y^\infty(\OO)$.

For, let $j,s\in\N$ and $(p,q)\in F$.

We first prove that the set $E(j,p,q,s)$ is open in $Y^\infty(\OO)$. Indeed, let $f\in E(j,p,q,s)$ and let $g\in Y^\infty(\OO)$ be such that $\dsup_{z\in L_m}|f^{(t)}(z)-g^{(t)}(z)|<a$ for $t=0,1,2,\ld,p+q+1$ (relation (1)).

The number $a>0$ will be determined later on. It is enough to prove that if $a$ is small enough then $g\in E(j,p,q,s)$.

The Hankel determinants defining $D_{p,q}(\zi)$ for $f$ depend
continuously on $\zi\in L_m$; thus, there exists $\de>0$ such that the
absolute values of the corresponding Hankel determinants are greater
than $\de>0$, for every $\zi\in L_m$, because $f\in D_{p,q}(\zi)$ for
$\zi\in L_m$ and because $L_m$ is compact.

From relation (1) we can control the first $p+q+1$ Taylor
coefficients of $g$ and by making $a>0$ small enough one can get the
Hankel determinants that define $D_{p,q}(\zi)$ to have absolute
value at least $\de/2>0$.

Therefore, $g$ will belong to $D_{p,q}(\zi)$ for every $\zi\in L_m$.
Now we consider the Pad\'{e} approximants of $f,g$ according to the
Jacobi formula (see preliminaries)
\[
[f;p/q]_\zi(z)=\frac{A(f,\zi)(z)}{B(f,\zi)(z)} \ \ \text{and} \ \ [g;p/q]_\zi(z)=\frac{A(g,\zi)(z)}{B(g,\zi)(z)}.
\]
Now $\mid A(f,\zi)(z)\mid^2+\mid B(f,\zi)(z)\mid^2$ vary
continuously with respect to $(z,\zi)\in K\times L_m$, because of the
Jacobi formula. So, there is a $\de'>0$ such that:
\[
\mid A(f,\zi)(z)\mid^2+\mid B(f,\zi)(z)\mid^2\ge\de', \ \ \text{for all} \ \ \zi\in L_m \ \ \text{and} \ \ z\in K.
\]
Now again from the Jacobi formula, if $a$ is small enough, one gets:
\[
\mid A(g,\zi)(z)\mid^2+\mid B(g,\zi)(z)\mid^2\ge\de'/2, \ \ \text{for all} \ \ \zi\in L_m \ \ \text{and} \ \ z\in K.
\]
This yields that $g\in E_{p,q,\zi}(K)$ for every $\zi\in L_m$. For the
rest it is enough to show that if $a$ is small enough then
$\dis\sup_{\zi\in L_m}\,\dis\sup_{z\in
K}\chi([g;p/q]_\zi(z),[f;p/q]_\zi(z))$ can become less than
$\dfrac{1}{s}-\dis\sup_{\zi\in L}\,\dis\sup_{z\in
K}\chi([f;p/q]_\zi(z),f_j(z))\equiv\ga>0$. By taking $a$ small as
before we have that $\mid A(f,\zi)(z)\mid^2+\mid
B(f,\zi)(z)\mid^2>\de'$ and for every $\zi\in L_m$ and $z\in K$ we
have $\mid A(g,\zi)(z)\mid^2+\mid B(g,\zi)(z)\mid^2>\de'/2$.

It follows that
\begin{align*}
\chi&([f;p/q]_\zi(x),[g;p/q]_\zi(z))\\
=&\,\frac{\mid A(f,\zi)(z)B(g,\zi)(z)-A(g,\zi)(z)B(f,\zi)(z)\mid}
{\sqrt{\mid A(f,\zi)(z)\mid^2+\mid B(f,\zi)(z)\mid^2}
\sqrt{\mid A(g,\zi)(z)\mid^2+\mid B(g,\zi)(z)\mid^2}} \\
&\le\frac{\sqrt{2}}{\de'}\mid A(f,\zi)(z)B(g,\zi)(z)-A(g,\zi)(z)B(f,\zi)(z)\mid
\end{align*}
for all $\zi\in L_m$ and $z\in K$, which easily yields the result, because the last expression can become as small as we want to, uniformly for all $\zi\in L_m$, $z\in K$.\ Thus, we proved that $E(j,p,q,s)$ is open.

Next, we prove that $T(p,q,s)$ is also open in $Y^\infty(\OO)$. Let $f$ be a function in $T(p,q,s)$ and let $g$ be a function in $Y^\infty(\OO)$ such that: $\dis\sup_{z\in L'_m}\mid f^{(t)}(z)-g^{(t)}(z)\mid<a$, for $t=0,1,2,\ld,\max(s,p+q+1)$, where $a>0$ will be determined later on.

In the same way as before one deduces that by making ``$a$'' small
enough it follows that $g\in D_{p,q}(\zi)\cap E_{p,q,\zi}(L_m)$, for all
$\zi\in L_m$.

Now $f(z)\in\C$, for each $z\in L_m$.\ It follows that for all $\zi\in L_m$, $z\in L_m$ we have $[f;p/q]_\zi(z)\in\C$. Therefore, $B(f,\zi)(z)\neq0$, where $B$ is given by the Jacobi formula.

So there is a $\de''>0$ such that $\de''<1$ and $\mid B(f,\zi)(z)\mid>\de''$ for all $\zi\in L_m$, $z\in L_m$, because $L_m\times L_m$ is compact.

By making ``$a$'' small enough, by continuity one can get
\[
\mid B(g,\zi)(z)\mid>\frac{\de''}{2} \ \ \text{for all} \ \ \zi\in L_m, \ \ z\in L_m.
\]
For $\el\in\{0,1,\ld,s\}$ it holds
\begin{align*}
\sup_{\zi\in L_m}\,\sup_{z\in L_m}\mid\big[g;p/q\big]_{\zi}^{(\el)}(z)-g^{(\el)}(z)
\mid\le&\,\sup_{z\in L_m}|f^{(\el)}(z)-g^{(\el)}(z)\mid \\
&+\sup_{\zi\in L_m}\,\sup_{z\in L_m}|f^{(\el)}(z)-\big[f;p/q\big]^{(\el)}_{\zi}(z)\mid\\
&+\sup_{\zi\in L_m}\,\sup_{z\in L_m}\mid\big[g;p/q\big]^{(\el)}_\zi(z)
-\big[f;p/q\big]^{(\el)}_\zi(z)\mid.
\end{align*}
The first term obviously get small as ``$a$'' gets small. Since the second term is fixed and less than $1/s$ we must control only the last term.

But the Jacobi denominators of $\big[f;p/q\big] ^{(\el)}_\zi(z)$ and
$\big[g;p/q\big]^{(\el)}_\zi(z)$ are bounded below from
$(\de'')^{\el+1}$ and $(\de''/2)^{\el+1}$ respectively for
$\el=0,1,\ld,s$.

Thus, the last term can get as small as we want for all $\el=0,1,\ld,s$, if $a$ is small enough.\ We are done.

Finally, we prove that for all $j,s\in\N$ the set $\dis\bigcup_{(p,q)\in F}(E(j,p,q,s)\cap T(p,q,s))$ is dense in $Y^\infty(\OO)$.

Let $g$ be a function inside $Y^\infty(\OO)$, $\e>0$ and $r,N\in\N$. Because of the definition of $Y^\infty(\OO)$ we can assume without loss of generality that $g$ is a rational function with poles off $\oO$, as we want to approximate $g$ on $L_r$ and every connected component of $\{\infty\}\cup(\C\sm L_r)$ contains a point from $\{\infty\}\cup(\C\sm\oO)$.

To prove what we want to, we have to find a function $f\in Y^\infty(\OO)$ and a pair $(p,q)\in F$ such that
\begin{enumerate}
\item[(1)] $f\in D_{p,q}(\zi)\cap E_{p,q,\zi}(K\cup L_m)$, for all $\zi\in L_m$.
\item[(2)] $\dsup_{\zi\in L_m}\dsup_{z\in K}\chi([f;p/q]_\zi(z),f_j(z))<\dfrac{1}{s}$.
\item[(3)] $\dsup_{\zi\in L_m}\dsup_{z\in L_m}\big|f^{(\el)}(z)-\big[f;p/q\big]^{(\el)}_\zi(z)\big|<\dfrac{1}{s}$, for $\el=0,1,\ld,s$.
\item[(4)] $\dsup_{z\in L_r}|f^{(t)}(z)-g^{(t)}(z)|<\e$, for $t=0,1,\ld,N$.
\end{enumerate}
Without loss of generality we may assume that $r>m$ or equivalently $L_r\supset L_m$.

Let $w:L_r\cup K\ra\C$ such that $q(z)=\left\{\begin{array}{cc}
                                                f_j(z) & z\in K \\
                                                g(z), & z\in L_r
                                              \end{array}\right.$.

Now, let $\mi$ be the sum of the principal parts of the poles of the
rational function $f_j$ that belong to $K$. Then $(\oo-\mi)$ is
holomorphic in a neighborhood of $(L'_r\cup K)$. Combining Runge's
with Weierstrass Theorems we conclude that there exists a rational
function $\dfrac{\tA(z)}{\tB(z)}$ with poles out of $(L'_r\cup K)$,
approximating $(\oo-\mi)$ uniformly on $L'_r\cup K$ with respect to
the euclidean metric and in the level of all derivatives of order
from zero to $N$. That implies that the function
$\dfrac{A(z)}{B(z)}=\mi(z)+\dfrac{\tA(z)}{\tB (z)}$ approximates
$f_j(z)$, uniformly on $K$ with respect to the chordal distance, and
also that $\Big(\dfrac{A(z)}{B(z)}\Big)^{(\el)}$ approximates the
function $(g(z))^{(\el)}$ uniformly on $L'_r$, with respect to the
euclidean metric.\ Obviously, we can assume that the greatest common
divisor of $A(z)$ and $B(z)$ is equal to one.

From our assumption on $F$, there exists a pair $(p,q)\in F$ such
that $p>\deg A,\deg B$ and $q>\deg B$. We consider the function
$\dfrac{A(z)}{B(z)}+dz^T=\dfrac{A(z)+d\cdot z^T\cdot B(z)}{B(z)}$
where $T=p-\deg B$ and $d$ is different from zero. Now, it is easy
to see that $gcd(A(z)+dz^TB(z),B(z))$ equals again to one.\ Thus,
according to Proposition \ref{prop2.2} it holds that for all
$\zi\in\C$ such that $B(\zi)\neq0$ the rational function
$\dfrac{A(z)+dz^TB(z)}{B(z)}$ belongs to $D_{p,q}(\zi)$ and also
$\Big[\dfrac{A(z)+dz^TB(z)}{B(z)};p/q\big]_\zi(z)=\dfrac{A(z)+dz^TB(z)}
{B(z)}$. In particular the above holds for all $\zi\in L'_m$, because
$B(\zi)\neq0$ for all $\zi\in L'_r$.

We distinguish the case $B(z)\neq0$ for all $z\in\dis\bigcup_nL_n$ and the case
where $B$ has zeros in $\Big(\dis\bigcup_{n\in\N}L_n\Big)\sm L_r$. First assume that $B(z)\neq0$
for all $z\in\dis\bigcup_nL_n$. In this case we set
$f(z)=\dfrac{A(z)}{B(z)}+dz^T$, and by selecting $d$ with $|d|$
small enough, we are done.

In the second case, since every component of $\widetilde{\C}\sm L_r$
contains a point from $\widetilde{\C}\sm\oO$, there exists a
rational function that belongs to $Y^\infty(\OO)$, call it $f$, such
that every finite set of derivatives $f^{(\el)}$ are close to
$\Big(\dfrac{A(z)+dz^T}{B(z)}\Big)^{(\el)}$ uniformly on $L_r$. This
follows immediately from Runge's and Weierstrass Theorems and also from the
fact that $B$ has finitely many roots outside $L_r$ and thus in a
positive distance from $L_r$.

It is easy to see that $f$ fulfills all requirements in the same way
as $\dfrac{A(z)+dz^TB(z)}{B(z)}$ does except from the fact that
maybe $[f;p/q]_\zi(z)\neq f(z)$.

But the following is true:
\begin{align*}
\dis\sup_{\zi\in L}\,\sup_{z\in\De}\mid\big[f;p/q\big]^{(\el)}_\zi(z)-f^{(\el)}(z)\mid\le&\,\sup_{z\in\De}
\mid f^{(\el)}(z)-h^{(\el)}(z)\mid \\
&+\sup_{\zi\in L}\,\sup_{z\in\De}\mid h^{(\el)}(z)-\big[h;p/q\big]^{(\el)}_\zi(z)\mid\\
&+\sup_{\zi\in L}\,\sup_{z\in\De}\mid\big[h;p/q\big]^{(\el)}_\zi(z)-
\big[f;p/q\big]^{(\el)}_\zi(z)\mid.   \hspace*{0.8cm} \mbox{($\ast$)}
\end{align*}
with $h(z)=\dfrac{A(z)+dz^TB(z)}{B(z)}$.

Now, as $p,q$ are fixed and we can control any finite set of
derivatives of $f$, we can also control any finite set of Taylor
coefficients of $f$. Thus, we can make the first and the last term
of the right-hand side expression in ($\ast$) as small as we want
and we are done.

This completes the proof of the Theorem. \qb
\end{Proof}

Now, if we define as above $\cu^K(m)$ the set of all the functions that satisfy the requirements of Theorem \ref{thm3.1} we have shown that for every $m\in\N$ and $K\subseteq\C$ compact such that $K\cap L_n=\emptyset$, for every $n\in\N$, $\cu^K(m)$ is a $G_\de$-dense subset of $Y^\infty(\OO)$. Therefore, if we have a sequence $(K_r)_{r\in\N}$ of such compact sets and apply Baire's Theorem at $\dis\bigcap_{r,m\in\N}\cu^{K_r}(m)$ we get the following result.
\begin{thm}\label{thm3.2}
Let $F\subset\N\times\N$ be a set that contains a sequence $(\tp_n,\tq_n)\in F$ such that $\tp_n\ra+\infty$, $\tq_n\ra+\infty$ and let $\OO\subseteq\C$ be an open set. Then there exists a function $f\in Y^\infty(\OO)$ such that for every rational function $h$ and $K_s$ member of the sequence $(K_r)_{r\in\N}$ there exists a sequence $(p_n,q_n)\in F$, $n=1,2,\ld$ with the following properties:

For every $m\in\N$ there exists $n(m)\in\N$ such that

(i) \ \  $f\in D_{p_n,q_n}(\zi)\cap E_{p_n,q_n\zi}(K_s\cup L_m)$, for every $n\ge n(m)$

(ii) For every $\el\in\N$, $\dsup_{\zi\in L_m}\dsup_{z\in L_M}\big|\big[f;p_n/q_n\big]^{(\el)}_{\zi}(z)-f^{(\el)}(z)\big|\ra0$, as $n\ra+\infty$ and

(iii) $\dsup_{\zi\in L_m}\dsup_{z\in K_s}\chi([f;p_n/q_n]_\zi(z)$, $h(z))\ra0$, as $n\ra+\infty$.

The set of such functions $f\in Y^\infty(\OO)$ is dense and $G_\de$ in $Y^\infty(\OO)$.
\end{thm}

Further, if we start with an open set $\OO$ such that $\{\infty\}\cup(\C\sm\oO)$ is connected then we can get the following result.
\begin{thm}\label{thm3.3}
Let $\OO$ be an open set such that $\{\infty\}\cup(\C\sm\oO)$ is connected and let $F$ be a subset of $\N\times\N$ that contains a sequence $(\tp_n,\tq_n)$, $n=1,2,\ld$ with $\tp_n\ra+\infty$.

Let also $K\subseteq\C$ be a compact set such that $K\cap L_n=\emptyset$ for every $n\in\N$ and let also $m\in\N$ be a fixed natural number.

Then there exists a function $f\in Y^\infty(\OO)$ such that for every polynomial $p$ there exists a sequence $(p_n,q_n)\in F$ $(n=1,2,\ld)$ with the following properties:\vspace*{-0.2cm}
\begin{enumerate}
\item[(i)] $f\in D_{p_n,q_n}(\zi)\cap E_{p_n,q_n,\zi}(L_m\cup K)$, $\forall\;\zi\in  L_m$.\vspace*{-0.2cm}
\item[(ii)] For every $\el\in\N$, $\dsup_{\zi\in L_m}\dsup_{z\in L_m}|[f;p_n/q_n]^{(\el)}(z)-f^{(\el)}(z)|\ra0$, as $n\ra+\infty$.\vspace*{-0.2cm}
\item[(iii)] For every $\el\in\N$, $\dsup_{\zi\in L_m}\dsup_{z\in K}\big|\big[f;p_m/q_n\big]^{(\el)}_\zi(z)-p^{(\el)}(z)\big|\ra0$, as $n\ra+\infty$.\vspace*{-0.2cm}
\end{enumerate}

The set of such functions $f\in Y^\infty(\OO)$ is dense and $G_\de$ in $Y^\infty(\OO)$.
\end{thm}
\begin{Proof}
First we enumerate the polynomials with coefficients from $\Q+i\Q$, $(\rho_j)$, $j=1,2,\ld$ and we call $V^K(m)$ the set of the functions that satisfy the prerequisites of the Theorem \ref{thm3.3}.

Now, we set for $j,s\in\N$ and $(p,q)\in\F$,
\begin{align*}
E(j,p,q,s)=\big\{&f\in Y^\infty(\OO)\mid f\in D_{p,q}(\zi)\cap E_{p,q,\zi}(K), \ \ \text{for every}\ \ \zi\in L_m  \\
&\text{and} \ \ \dsup_{\zi\in L_m}\dsup_{z\in K}|[f_ip/q]_\zi(z)-p_j(z)|<\frac{1}{s}\big\}
\end{align*}
and
\begin{align*}
T(p,q,s)=\big\{&f\in Y^\infty(\OO)\mid f\in D_{p,q}(\zi)\cap E_{p,q,\zi}(L_m), \ \ \text{for every} \ \ \zi\in L_m \\
&\text{and} \ \ \dsup_{\zi\in L_m}\dsup_{z\in L_m}\big|\big[f;p/q\big]^{(\el)}_\zi(z)-f^{\el)}(z)\big|<\frac{1}{s} \ \ \text{for} \ \ \el=0,1,\ld,s\big\}.
\end{align*}
Then the definition of $V^K(m)$ yields that
\[
V^K(m)=\bigcap^{+\infty}_{j,s=1}\bigcup_{(p,q)\in F}(E(j,p,q,s)\cap T(p,q,s)).
\]
Then, we prove that $V^K(m)$ is a $G_\de$-dense subset of $Y^\infty(\OO)$ and therefore $V^K(m)\neq\emptyset$, which is what we want to.

For, we prove that for every $i,s=1,2,\ld$ and $(p,q)\in F$ the sets $E(j,p,q,s)$, $T(p,q,s)$ are open in $Y^\infty(\OO)$ and that for every $j,s=1,2$ the set $\dis\bigcup_{(p,q)\in F}(E(j,p,q,s)\cap T(p,q,s))$ is dense in $Y^\infty(\OO)$.

The proof that $E(j,p,q,s)$ and $T(p,q,s)$ are open in $Y^\infty(\OO)$ is the same as in the corresponding of Theorem \ref{thm3.1} and therefore is omitted.

To prove the density result let $i,s=1,2,\ld\;.$

We will show that $\dis\bigcup_{(p,q)\in F}(E(j,p,q,s)\cap T(p,q,s)$ is dense in $Y^\infty(\OO)$.

For, let $g$ be a rational function with poles off $\oO$, $\e>0$, $N\in\N$ and $L_r$ be a member of the sequence of compacts $(L_n)_{n\in\N}$ that defined $Y^\infty(\OO)$. Because $(L_n)_{n\in\N}$ is an increasing sequence of compact sets we may assume that $L_r\supset L_m$ and because of the fact that $\{\infty\}\cup(\C\sm\oO)$ is connected and Runge's Theorem we may also assume that $q$ is a polynomial.

Now to prove what we want to, we need to find a function $f\in Y^\infty(\OO)$ and a pair $(p,q)\in F$ such that\vspace*{-0.2cm}
\begin{enumerate}
\item[(i)] $f\in D_{p,q}(\zi)\cap E_{p,q,\zi}(K\cup L_m), \ \ \forall\;\zi\in L_m$.\vspace*{-0.2cm}
\item[(ii)] $\dsup_{\zi\in L_m}\dsup_{z\in K}|[f;p/q]_\zi(z)-f_j(z)|<\dfrac{1}{s}$.\vspace*{-0.2cm}
\item[(iii)] $\dsup_{\zi\in L_m}\dsup_{z\in L_m}\big|\big[f;p/q\big]^{(\el)}_\zi(z)-f^{(\el)}(z)\big|<\dfrac{1}{s}$, for $\el=0,1,\ld,s$ and\vspace*{-0.2cm}
\item[(iv)] $\dsup_{z\in L_r}|f^{(t)}(z)-g^{(t)}(z)|<\e$, for $t=0,1,\ld,N$.\vspace*{-0.2cm}
\end{enumerate}

Let $w:L_r\cup K\ra\C$ defined by $w(z)=\left\{\begin{array}{cc}
                                                 f;(z), & z\in K \\
                                                 g(z), & z\in L_r
                                               \end{array}\right.$. Because of the fact that every connected
component of $\{\infty\}\cup(\C\sm L_r)$ contains a connected component of $\{\infty\}\cup(\C\sm\oO)$ and the last one is a connected set we have that $\{\infty\}\cup(\C\sm L_r)$ is also connected.

Moreover, since $K^c$ is connected and $K\cap L_n=\emptyset$ is a compact set Runges and Weierstrass theorem yields the existence of a polynomial $h$ that  approximates
 $w$ in the level of any finite set of derivatives on $K\cup L_r$.
Now, there exist a pair $(p,q)\in F$ with $p>\deg(g)$.

Therefore, if we set $f(z)=g(z)+dz^p$, for $d\in\C$ with $|d|\neq0$ small enough one can check with the help of Proposition \ref{prop2.2} that $f(z)$ satisfies all the four conditions we want. This completes the proof of Theorem 3. \qb
\end{Proof}

Now, if we set $V^K(m)$, as in the above proof, for $K\subseteq\C$ a compact set with connected complement such that $K\cap L_n=\emptyset$, for every $n\in\N$ and $m\in\N$, the set of all functions $f\in Y^\infty(\OO)$ that satisfy the requirements of Theorem 3, we have shown that $V^K(m)$ is a $G_\de$ and dense subset of $Y^\infty(\OO)$. Therefore, by applying Baire's Theorem at $\dis\bigcap_{\el,m\in\N}V^{K_\el}(m)$, where $(K_\el)_{\el\in\N}$ is a set of compact sets in $\C$ with connected complement such that $K_\el\cap L_n=\emptyset$ for every $\el,n\in\N$ we get the following result.
\begin{thm}\label{thm3.4}
Let $F\subset\N\times\N$ be a set that contains a sequence $(\tp_n,\tq_n)\in F$ such that $\tp_n\ra+\infty$ and let $\OO\subseteq\C$ be an open set. Then there exists a function $f\in Y^\infty(\OO)$ such that for every polynomial function $p$ and $K_s$ member of the sequence $(K_\el)_{\el\in\N}$ that we defined above, there exist a sequence $(p_n,q_n)\in F$, $n=1,2,\ld$ with the following properties: For every $m\in\N$, there exists $n(m)\in\N$ such that:\vspace*{-0.2cm}
\begin{enumerate}
\item[(i)] $f\in D_{p_n,q_n}(\zi)\cap E_{p_n,q_n,\zi}(K_s\cup L_m)$, for every $n\ge n(m)$.\vspace*{-0.2cm}
\item[(ii)] For every $\el\in\N$, $\dsup_{\zi\in L_m}\dsup_{z\in L_m}\big|\big[f;p_n/q_n\big]^{(\el)}_\zi(z)-f^{(\el)}(z)\big|\ra0$.\vspace*{-0.2cm}
\item[(iii)] $\dsup_{\zi\in L_m}\dsup_{z\in K_s}|[f;p_n/q_n]_\zi(z)-p(z)|\ra0$.\vspace*{-0.2cm}
\end{enumerate}

The set of such functions $f\in Y^\infty(\OO)$ is dense and $G_\de$ in $Y^\infty(\OO)$.
\end{thm}
\begin{rem}\label{rem3.5}
One can check that in the case where $\OO$ is a simply connected domain and $L_n\cap\partial\OO=\emptyset$, for every $n\in\N$, the proofs of Theorem \ref{thm3.3}, \ref{thm3.4} can be easily modified to show that the results under this new hypothesis remain true.
\end{rem}
\section{Applications}\label{sec4}
\noindent

We consider now different cases for the sequences $(L_n)_{n\in\N}$ and $(K_r)_{r\in\N}$ and we discuss what old or new results one can get from the theorems of Section \ref{sec3} under these hypothesis.\vspace*{0.2cm} \\
\noindent
{\bf 4.1.} Consider first the case where $(L_n)_{n\in\N}$ is an exhausting sequence of compact sets in $\OO$, where $\OO\subseteq\C$ is an open set and $K=\emptyset$.

Such a sequence $(L_n)_{n\in\N}$ exists because of the fact that $\OO$ is an open set and of Proposition \ref{prop2.7}. Moreover, such a sequence can be easily modified as in Corollary \ref{cor2.9} to give also an exhausting sequence of compact sets in $\OO$ and to fulfill also the conditions to define the space $T^\infty(\OO)=T^\infty(\OO,(L_n)_{n\in\N})$, that we introduced at Section \ref{sec3}. Furthermore, in this case $T^\infty(\OO)$ becomes equal to $H(\OO)$. But by Runge's Theorem the set of rational function with poles off $\OO=\dis\bigcup_n L_n$ is dense in $H(\OO)$ which can be translated to $Y^\infty(\OO)=H(\OO)$, where $Y^\infty(\OO)$ is also defined in Section \ref{sec3}{. Now, under these hypothesis, Theorem \ref{thm3.1} yields a result from \cite{6}, according to which generically every holomorphic function $f$ defined on an open set $\OO$ can be approximated by a sequence of it's Pad\'{e} approximants obtained from a set $\{[f;p/q]\mid(p,q)\in F\}$ where $F$ is a subset of $(\N\times\N$ that contains a sequence $(\tp_n,\tq_n)$, $n=1,2,\ld$ where $\tp_n\ra+\infty$, $\tq_n\ra+\infty$.\vspace*{0.2cm} \\
\noindent
{\bf 4.2.} We consider now the case where $(L_n)_{n\in\N}$ is an exhausting sequence of compact sets in $\OO$, where $\OO\subseteq\C$ is an open set, and
\[
K_r=(\C\sm\OO)\cap\{z\in\C\mid|z|\le r\}, \ \ r\in\N.
\]
Then, as in Application 4.1, $(L_n)_{n\in\N}$ defines the space $T^\infty(\OO,L_n)=T^\infty(\OO)$ which, as a set, equals $H(\OO)$ and furthermore, it holds $Y^\infty(\OO)=T^\infty(\OO)=H(\OO)$. Under these hypothesis Theorem \ref{thm3.2} yields Theorem \ref{thm3.4} of \cite{13}, according to which if $F$ is a subset of $\N\times\N$ that contains a sequence $(\tp_n,\tq_n)$, $n=1,2,\ld$ with $\tp_n\!\!\ra\!\!+\infty$, $q_n\!\!\ra\!\!+\infty$, then generically every function $f\in H(\OO)$ has universal Pad\'{e} approximants obtained from the set $\{[f;p/q]\mid(p,q)\in F\}$ that can approximate uniformly any rational function on any compact subset of $\C\sm\OO$ with respect to the chordal distance and can, simultaneously, approximate uniformly $f$ on any compact subset of $\OO$ with respect to the Euclidean distance. One can see that in this case the universal approximation is valid on $\partial\OO$.\vspace*{0.2cm} \\
\noindent
{\bf 4.3.} We consider now the case where $\OO\subseteq\C$ is a simply connected domain with a locally finite number of components, $(L_n)_{n\in\N}$ is an exhausting family of compact sets in $\OO$ (Proposition \ref{prop2.7}) and $(K_r)_{r\in\N}$ is a sequence of compact sets in $\C\sm\OO$ with connected complement as in Proposition \ref{prop2.3}. The sequence $(L_n)_{n\in\N}$, as in Application 4.1, defined the space $T^\infty(\OO)=T^\infty(\OO,(L_n)_{n\in\N})$ and it holds $T^\infty(\OO)=Y^\infty(\OO)=H(\OO)$. Now, under these hypothesis, if the set $F\subset\N\times\N$ contains a sequence $(\tp_n,\tq_n)$ with $\tp_n\ra+\infty$, then Theorem \ref{thm3.4} combined with Remark \ref{rem3.5} yields Theorem \ref{thm3.4} of \cite{5}, according to which generically every holomorphic function $f\in H(\OO)$ has universal Pad\'{e} approximants obtained from the set $\{[f;p/q]\mid(p,q)\in F\}$ that can approximate any polynomial function on any compact subset of $\C\sm\OO$ with connected complement with respect to the Euclidean distance and can, simultaneously, approximate $f$ with respect to the Euclidean distance on every compact subset of $\OO$. We see that the universal approximation is also in this case valid on $\partial\OO$.\vspace*{0.2cm}
\noindent
{\bf 4.4.} We consider now the case where $L_n=\oO\cap\overline{D(0,n)}$ for every $n\in\N$ and $(K_r)_{r\in\N}$ is an exhausting sequence of compact sets in $\C\sm\oO$. The existence of the sequence $(K_r)_{r\in\N}$ is established from Proposition \ref{prop2.7}.

Now, it is easy to see that $(L_n)_{n\in\N}$ is an increasing sequence of compact subsets of $\oO$ such that it holds $\overline{(L_n\cap\OO)}=L_n$, for every $n\in\N$ and for every compact set $J\subset\OO$, there exists $n\in\N$ such that $J\subset\oO\cap\overline{D(0,n)}=L_n$. For the last condition, let $n\in\N$. We will show that every connected component of $\{\infty\}\cup(\C\sm L_n)$ contains a point from $\{\infty\}\cup(\C\sm\oO)$. It holds that $\{\infty\}\cup(\C\sm L_n)=\{\infty\}\cup(\C\sm\oO)\cup(\C\sm\overline{D(0,n)})$, and if $C$ is a connected component of $\{\infty\}\cup(\C\sm L_n)$ it holds either $C\subset\{\infty\}\cup(\C\sm\oO)$ or $C\cap(\C\sm\overline{D(0,n)})\neq\emptyset$. In the former case the component $C$ obviously contains a point from $\{\infty\}\cup(\C\sm\oO)$ and in the latter case $C$ must contain the whole connected set $\{\infty\}\cup(\C\sm\overline{D(0,n)})$ which means that $\{\infty\}\in C$ which is also a point from $\{\infty\}\cup(\C\sm\oO)$. Therefore, the space $T^\infty(\OO)\!\!=\!\!T^\infty(\OO,(L_n)_{n\in\N})$ can be defined and because $\dis\bigcup_nL_n$ contains any compact subset of $\oO$, $T^\infty(\OO)$ becomes equal to the space $A^\infty(\OO)$, which is the space of the functions $f\in H(\OO)$ such that for every $\el=0,1,\ld$ $f^{(\el)}$ extends continuously on $\oO$, and $Y^\infty(\OO)$ becomes equal to $X^\infty(\OO)$, which is the closure of all rational functions with poles off $\oO$. Under these hypothesis Theorem \ref{thm3.2} yields Theorem \ref{thm3.3} of \cite{14}, according to which if $F$ is a subset of $\N\times\N$ that contains a sequence $(\tp_n,\tq_n)$ where $\tp_n\!\!\ra\!\!+\infty$, $\tq_n\!\!\ra\!\!+\infty$ then generically every function $f\in X^\infty(\OO)\subseteq A^\infty(\OO)$ has universal Pad\'{e} approximants obtained from the set $\{[f;p/q]\mid(p,q)\in F\}$ that can approximate uniformly any rational function on any compact subset of $\C\sm\oO$ with respect to the chordal metric, and can, simultaneously, approximate $f$ on every compact subset of $\oO$ with respect to the Euclidean distance.

We see that in this case, not only the approximation is not valid on the boundary $\partial\oO$ but on the compact subsets of $\partial\OO$ the universal Pad\'{e} approximants approximate the function that defines them.\vspace*{0.2cm}\\
\noindent
{\bf 4.5.} Consider now the case where $\OO\subseteq\C$ is an open set with a locally finite number of components such that $\{\infty\}\cup(\C\sm\oO)$ is connected, $L_n=\oO\cap\overline{D(0,n)}$ for every $n\in\N$ and $(K_r)_{r\in\N}$ is an exhausting sequence of compact subsets of $\C\sm\oO$ with connected complement as in Proposition \ref{prop2.4}.

Then, as in Application 4.4, $T^\infty(\OO)=T^\infty(\OO,(L_n)_{n\in\N})$ coincides with the space $A^\infty(\OO)$ and $Y^\infty(\OO)$ coincides with $X^\infty(\OO)$, where both $A^\infty(\OO)$ and $X^\infty(\OO)$ were defined in Application 4.4. Then Theorem \ref{thm3.3} under these hypothesis yields that if $F$ is a subset of $\N\times\N$ that contains a sequence $(\tp_n,\tq_n)$, $n=1,2,\ld$ with $\tp_n\ra+\infty$, then generically every function $f\in X^\infty(\OO)\subseteq A^\infty(\OO)$ has universal Pad\'{e} approximants obtained from the set $\{[f;p/q]\mid(p,q)\in F\}$ that can approximate uniformly every polynomial function on every compact subset of $\C\sm\oO$ with connected complement with respect to the Euclidean metric and can simultaneously approximate $f$ on every compact subset of $(\oO)$. The above result is similar to the Theorem \ref{thm4.1} of \cite{14}.

One can notice that also in this case the universal approximation is not valid on $\partial\OO$.\vspace*{0.2cm} \\
\noindent
{bf 4.6.} Now, we will present an example where on a part of the boundary $\partial\OO$ the universal Pad\'{e} approximation will be valid and on another part it will be not.
\[
\text{Let} \ \ \OO=D=\{z\in\C\mid|z|<1\}, \ \ \text{and for every} \ \ n\in\N \ \ \text{we set}
\]
\[
L_n=\{z\in\C\mid|z|\le1-1/n\}\cup\{z\in\C\sm\{0\}\mid|z|\le1,\;1/n\le\Arg z\le\pi-1/n\}.
\]
We also consider $(K_r)$ $r=1,\ld$ a sequence of compact subsets of $\C\sm\oD$ as in Proposition \ref{prop2.4}, and
\[
K_1=\{z\in\C\mid|z|\ge1,\;-2\le\Re(z)\le2,\;\text{Im}(z)\le0\}.
\]
It is easy to check that with this choice of $(L_n)_{n\in\N}$ the space $T^\infty(\OO)=T^\infty(\OO,(L_n)_{n\in\N})$ is well defined and that for this space Theorem \ref{thm3.4} yields the following theorem.
\begin{thm}\label{thm4.1}
Let $F$ be a subset of $\N\times\N$ such that $F\supset\{(\tp_n,\tq_n)$, $n\in\N$ with $\tp_n\ra+\infty\}$ and let $(L_n)_{n\in\N}$ and $(K_r)_{r\in\N}$ be defined as above. Then, there exist a function $f\in H(D)$ such that for every polynomial $p$ and $r\in\N$ there exists a sequence $(p_n,q_n)\in F$, $n=1,2,\ld$ such that for every member of the sequence $(L_n)_{n\in\N}$, $L_N$ there exists a natural number $m=m(N)\in\N$ such that:\vspace*{-0.2cm}
\begin{enumerate}
\item[(i)] $f\in D_{p_n,q_n}(\zi)\cap E_{p_n,q_n,\zi}(K_r\cup L_N)$ for every $n\ge m$.\vspace*{-0.2cm}
\item[(ii)] $\dsup_{\zi\in L_N}\dsup_{z\in K_r}|[f;p_n/q_n]_\zi(z)-p(z)|\ra0$, as $n\ra+\infty$.\vspace*{-0.2cm}
\item[(iii)] For every $\el\in\N_0$, $\dsup_{\zi\in L_N}\dsup_{z\in L_N}\big|\big[f;p_n/q_n\big]^{(\el)}_\zi(z)-f^{(\el)}(z)\big|\ra0$, as $n\ra+\infty$.\vspace*{-0.2cm}
\end{enumerate}
\end{thm}

One can see that Theorem \ref{thm4.1} in case $F=\{(n,0)\mid n\in\N\}$ coincides with Theorem \ref{thm3.1} of \cite{8}.

Now, Theorem \ref{thm4.1} for $r=1$, gives the existence of a function $f\in H(D)$ with Pad\'{e}\linebreak approximant $[f;p/q]$, $(p,q)\!\in\! F$, where $F\!\subset\!\N\!\times\!\N$ was defined above, that can approximate any polynomial function uniformly on any compact subset of $\{e^{i\vthi}\mid0<\vthi<\pi\}$\linebreak and can, simultaneously, approximate $f$ uniformly on any compact subset of\linebreak $\{e^{i\vthi}\mid\pi\!\le\!\vthi\!<\!2\pi\}\cup\{1\}$.

Therefore indeed we prove the existence of a function where the universal Pad\'{e} approximation is valid on a part of the boundary $\partial\OO$, for $\OO=D$, and it is not valid on another part of it.
\section{A sufficient condition}\label{sec5}
\noindent

Considering Theorem \ref{thm4.1} one can naturally ask whether one can find necessary and sufficient conditions about the subsets $S,T$ of the boundary $\partial\OO$, of an open set $\OO$, so that some Pad\'{e} approximants of a holomorphic function $f\in H(\OO)$ approximate universally any rational function on any compact subset of $(\C\sm\oO)\cup T$ with respect to the chordal distance and can, simultaneously, approximate $f$ on any compact subset of $\OO\cup S$ with respect to the Euclidean distance. For convenience, we denote the set of such functions by $V(S,T,\OO)$, possibly empty.

A necessary condition for $V(S,T,\OO)\neq\emptyset$ is of course $S\cap T=\emptyset$. We will prove now that a sufficient condition for $S,T\subseteq\partial\OO$ so that $V(S,T,\OO)\neq\emptyset$ is $\oS\cap\oT=\emptyset$; i.e. $S,T$ have disjoint closures.
\begin{prop}\label{prop5.1}
Let $\OO\subseteq\C$ be an open set and $S,T\subseteq\partial\OO$ such that $\oS\cap\oT=\emptyset$. Then $V(S,T,\OO)\neq\emptyset$.
\end{prop}
\begin{Proof}
Based on Theorem \ref{thm3.2}, it is enough to find two sequences $(L_n)_{n\in\N}$, $(K_r)_{r\in\N}$ of compact sets in $\oO$ and in $\C\sm\OO$, respectively, such that: Firstly, for every $r,n\in\N$ it holds $K_r\cap L_n=\emptyset$, secondly $(L_n)_{n\in\N}$ fulfills the prerequisites to define the space $T^\infty(\OO)=T^\infty(\OO,(L_n)_{n\in\N})$ and finally for any choice of compact subsets $J_1\subset\OO\cup S$, $J_2\subset(\C\sm\oO)\cup T$ there exist $R,N\in\N$ such that $J_1\subset L_N$ and $J_2\subset K_R$.

Because $\oS,\oT$ are closed, and therefore locally compact subsets of $\partial\OO$, from Proposition \ref{prop2.7} there exists an exhausting sequence of compact sets in $\oS$, call it $(L'_n)_{n\in\N}$ and an exhausting sequence of compact sets in $\oT$, call it $(K'_r)_{r\in\N}$. Now, as $\oS\cap\oT=\emptyset$, for every $n\in\N$ there exists a positive number $a_n>0$, such that $2a_n<\min(\dist(L'_n\oT)$, $\dist(\oS,K'_n))$.

Let also $(L''_n)_{n\in\N}$, $(K''_r)_{r\in\N}$ be two exhausting sequences of compact sets in $\OO$ and $\C\sm\oO$ respectively; the existence of them is established from the fact that they are open subsets of $\C$ combined with Proposition \ref{prop2.7}. Now, we set for every $n\in\N$,
\[
\left.\begin{array}{l}
L_n=\bigg(\bigcup^n_{N=1}\{z\in\oO\mid\dist(z,L'_N)\le a_N\}\bigg)\cup L''_n   \\
\text{and for every}\ \  r\in\N, \\
  K_r=\bigg(\bigcup^r_{R=1}\{z\in\oO\mid\dist(z,K'_R)\le a_R\bigg)\cup K''_r.
\end{array}\right\} \ \ \text{Definitions} \; (\ast)
\]
We first prove that for every $r,n\in\N$, $K_r\cap L_n=\emptyset$. Indeed, let $z\in L_n\cap K_r$. Then obviously $z\in\partial\OO$ and as $L''_n\subset\OO$ and $K''_r\subset\C\sm\oO$ we get that there exist $N\in\{1,2,\ld,n\}$ and $R\in\{1,2,\ld,r\}$ such that
\[
z\in\{z\in\oO\mid\dist(z,L'_N)\le a_N,\;\dist(z,K'_R)\le a_R\}.
\]
This means that there exist $\el_N\in L_N$, $k_R\in K_R$ such that $|z-\el_N|\le a_N$, $|z-k_R|\le a_R$ which gives $|\el_N-k_R|\le a_R+a_N$ and therefore $|\el_N-k_R|<2\max(a_R,a_N)$.

Without loss of generality, we may assume $a_N>a_R$ and therefore that it holds $|\el_N-k_R|\le2a_N$. The above relation yields $\dist(L'_N,\oT)\le\dist(L'_N,K'_R)\le2a_N$, which is a contradiction. Therefore, indeed for every $n,r\in\N$, $L_n\cap K_r=\emptyset$.

Now, it is easy to see that $(L_n)_{n\in\N}$ is an increasing sequence of compact subsets of $\oO$, that for every $n\in\N$ it holds $\overline{(L_n\cap\OO)}=L_n$, and that for every compact subset of $\OO$, $J\subset\OO$, there exists $N\in\N$ such that $J\subset L''_N\subset L_N$.

These three facts and Corollary \ref{cor2.9} show that $(L_n)_{n\in\N}$ can be modified such that for every $n\in\N$ every connected component of $\{\infty\}\cup(\C\sm L_n)$ will contain a point from $\{\infty\}\cup(\C\sm\oO)$, leaving however, for every $n\in\N$ the set $L_n\cap\partial\OO$ unchanged.

One can check that the new sequence $(L_n)_{n\in\N}$ fulfills all the prerequisites to define the space $T^\infty(\OO)=T^\infty(\OO,(L_n)_{n\in\N})$ and that it remains true that $L_n\cap K_r=\emptyset$, for every $r,n\in\N$. The former is immediate. For the latter; let $r,n\in\N$. Then, as $L_n\subset\oO$, $K_r\subset\C\sm\OO$ it must also be true that $L_n\cap K_r\subset\partial\OO$ and therefore $L_n\cap K_r=(L_n\cap\partial\OO)\cap(K_r\cap\partial\OO)$. The last equation together with the fact that $(L_n\cap\partial\OO)$ remained unchanged through the modification yields that indeed the intersection $L_n\cap K_r$ also remains equal to the empty set. Finally, we have to show that for any choice of compact subsets $J_1\subset\OO\cup S$ and $J_2\subset(\C\sm\oO)\cup T$, there exist $R,N\in\N$ such that $J_1\subset L_N$, $J_2\subset K_R$. As the new sequence $(L_n)_{n\in\N}$ contains the old one, we may use for $(L_n)_{n\in\N}$ the one we defined at $(\ast)$.

Now, because of the symmetry of the definitions $(\ast)$, it is enough to prove that for every compact set $J_1\subset\OO\cup S$, one can find $N\in\N$ such that $J_1\subset L_N$. Indeed, let $J_1$ be a compact subset of $\OO\cup S$. One can see that because $J_1\cap S$ is a compact subset of $\partial\OO$ and $(L'_n)_{n\in\N}$ is an exhausting sequence of compact sets in $\partial\OO$, there exists $n_1\in\N$ such that $J_1\cap S\subset L_{n_1}$. (Proposition \ref{prop2.6}).

Now, we claim that there exists $M\in\N$ such that
\[
J_1\cap\{z\in\oO\mid\dist(z,\partial\OO)<1/M\}\subseteq\{z\in\oO\mid\dist(z,L'_{n_1})\le a_{n_1}\}.
\]
Indeed, if not there exists a sequence $(x_n)_{n\in\N}\subset J_1$ such that $\dist(x_n,\partial\OO)\ra0$ and $\dist(x_n,L'_{n_1})\ge a_{n_1}$, for every $n\in\N$. These two relations together with the fact that $J_1$ is compact yield a point $z_0\in J\cap\partial\OO$ with $\dist(z_0,L'_{n_1})\ge a_{n_1}$. But this is a contradiction since $J_1\cap\partial\OO=J_1\cap S\subset L'_{n_1}$.

Now, fix such a number $M\in\N$. We have that $J_1\cap\{z\in\oO\mid\dist(z,\partial\OO)<1/M\}\subset L_{n_1}$, by definitions of $M$, $L_{n_1}$. Moreover, since $J_1\cap\{z\in\oO\mid\dist(z,\partial\OO)\ge1/M\}$ is a compact subset of $\OO$, there exists $n_2\in\N$ such that
\[
J_1\cap\{z\in\oO\mid\dist(z,\partial\OO)\ge1/M\}\subset L''_{n_2}\subset L_{n_2},
\]
because $(L''_n)_{n\in\N}$ is an exhausting sequence of compact sets in $\OO$.

Therefore, it holds:
\begin{align*}
J_1=&(J_1\cap\{z\in\oO\mid\dist(z,\partial\OO)<1/M)\cup(J_1\cap\{z\in\oO\mid\dist(z,\partial\OO)\ge1/M)\\
&\subseteq L_{n_1}\cup L_{n_2}\subset L_{n_1+n_2}=L_N,
\end{align*}
for $N=n_1+n_2\in\N$, which is what we wanted. The proof of Proposition \ref{prop5.1} is complete.
\qb
\end{Proof}
\begin{rem}\label{rem5.2}
In the case where $\OO\subseteq\C$ is a bounded open set, then one can simplify the above proof of Proposition \ref{prop5.1}, be replacing each $L'_n$ by $\oS$ and each $K'_r$ by $\oT$.
\end{rem}
\begin{cor}\label{cor5.3}
Let $\OO$ be an open set and $S,T\subseteq\C$ be two closed subsets of $\partial\OO$. Then, it holds $V(S,T,\OO)\neq\emptyset$ if and only if $S\cap T=\emptyset$.
\end{cor}
\begin{cor}\label{cor5.4}
Let $\OO$ be an open set. Then, if $S\subseteq\partial\OO$ is a closed set, there exists a function $f\in H(\OO)$, such that for every $\el\in\N$, the derivative $f^{(\el)}$ extends continuously on $\OO\cup S$.
\end{cor}

The total characterization of the subsets $S,T$ of the boundary $\partial\OO$ of an open set $\OO$ such that $V(S,T,\OO)\neq\emptyset$ remains an open problem.\vspace*{0.2cm} \\
{\bf Acknowledgement}. I would like to thank Professor V. Nestoridis for many helpful discussions and guidance through the creation of this paper.

\vspace*{1cm}
University of Athens \vspace*{0.1cm} \\
Department of Mathematics \vspace*{0.1cm} \\
157 84 Panepistimiopolis \vspace*{0.1cm} \\
Athens \vspace*{0.1cm} \\
GREECE \vspace*{0.1cm} \\
email address:ilias\_\_91@hotmail.com


\begin{thebibliography}{99}
\bibitem{1} G. A. Baker, Jr and P. R. Graves-Morris, Pad\'{e} Approximants. Bol. 1 and 2 (Encyclopedia of Math. and Applications), Cambridge, Un. Press 2010.
%
\bibitem{2} Bayart, Grosse-Erdmann, Nestoridis and Papadimitropoulos, Abstract theory of Universal series and applications, Proceedings of the Lon Math. Soc. (3) 96(2008) no 2, 417-463.
%
\bibitem{3} C. Chui and M. N. Parnes, Approximation by overconvergence of power series, J. Math. Anal. Appl. 36(1971), 693-696.
%
\bibitem{4} N. Daras, V. Nestoridis and Ch. Papadimitropoulos, Universal Pad\'{e} approximation of Selezner type, Arch. Math. (Basel), 100(2013) no 6 571-585.
%
\bibitem{5} N. Daras, G. Fournodavlos, V. Nestoridis, Universal Pad\'{e} approximants on simply connected domains, submitted.
%
\bibitem{6} G. Fournodavlos, V. Nestoridis, Generic approximation of functioning by their Pad\'{e} approximants, Journal of Mathematical Analysis and Application, vol. 408, issue 2 (2013), 744-750.
%
\bibitem{7}Ch. Kariofillis, Ch. Konstadilaki and V. Nestoridis, Smooth universal Taylor Series, Monasth. Math. 147(2006) no. 3, 249-257.
%
\bibitem{8} Ch. Kariofillis, V. Nestoridis: Universal Taylor series in Simly Connected Domains, C.M.F.T., Volume 6 (2006), No. 2, 437-446.
%
\bibitem{9} W. Luh (1970), Approximation analytisher Funktionen d\"{u}rch \"{u}berkonvergente Potenzreihen und dere Matrix-Transformierten. Mitt Math Sem Giessen 88:1-56.
%
\bibitem{10} A. Melas and V. Nestoridis: Universality of Taylor series as a generic property of holomorphic functions, Adv. Math. Vol. 157, 2001, p. 138-176.
%
\bibitem{11} V. Nestoridis, Universal Taylor series, Ann. Inst. Fourier 46 5(1996), 1293-1306.
%
\bibitem{12} V. Nestoridis, An extension of the notion of universal Taylor series in: N. Papamichael, S. Ruscheweyh, E. B. Saff (eds), Proceedings of the 3rd CMFT conference on computational methods and function theory 1997, Nicosia, Cyprus, October 13-17, 1997, World Scientific. Ser. Approx. Decompos. 11(1999), 421-430.
%
\bibitem{13} V. Nestoridis, Universal Pad\'{e} Approximants with respect to the chordal metric, Journal of Contemporary Mathematics Analysis, v. 47(2012) no 4, 168-181.
%
\bibitem{14} V. Nestoridis, I. Zadik, Pad\'{e} approximants, density of rational functions in $A^\infty(\OO)$ and smoothness of the integration operator, arxiv:1212.4394.
%
\bibitem{15} J. P\'{a}l (1914-1915) Zwei Kleine Bemerkungen, T\^{o}hoku Math J 6: 42-43.
\end{thebibliography}
\end{document}